\documentclass[aop,citesort,MSNbibl,dvips]{arximspdf}

\doi{10.1214/10-AOP527}
\volume{38}
\issue{4}
\pubyear{2010}
\firstpage{1672}
\lastpage{1689}

\makeatletter

\newtheorem{theorem}{Theorem}[section]
\newtheorem{proposition}{Proposition}[section]
\newtheorem{lemma}{Lemma}[section]
\newproclaim{remark}{Remark}[section]

\newcommand{\dd}{\stackrel{d}{\Rightarrow}}
\newcommand{\as}{\stackrel{\mathrm{a.s.}}{\rightarrow}}

\makeatother

\begin{document}
\begin{frontmatter}

\title{Bounds on the constant in the mean central limit~theorem}
\runtitle{On the mean central limit theorem constant}

\begin{aug}
\author[A]{\fnms{Larry} \snm{Goldstein}\corref{}\ead[label=e1]{larry@math.usc.edu}}
\runauthor{L. Goldstein}
\affiliation{University of Southern California}
\address[A]{Department of Mathematics, KAP 108\\
University of Southern California\\
Los Angeles, California 90089-2532\\
USA\\
\printead{e1}} 
\end{aug}

\received{\smonth{11} \syear{2008}}
\revised{\smonth{1} \syear{2010}}

%
\begin{abstract}
Let $X_1,\ldots,X_n$ be independent with zero means, finite variances
$\sigma_1^2,\ldots,\sigma_n^2$ and finite absolute third moments. Let
$F_n$ be the distribution function of
$(X_1+\cdots+X_n)/\sigma$, where $\sigma^2=\sum_{i=1}^n \sigma_i^2$,
and $\Phi$ that of the
standard normal. The $L^1$-distance between $F_n$ and
$\Phi$ then satisfies
\[
\Vert F_n-\Phi\Vert_1 \le\frac{1}{\sigma^3}\sum_{i=1}^n E|X_i|^3.
\]
In particular, when $X_1,\ldots,X_n$ are identically distributed with
variance $\sigma^2$,
we have
\[
\Vert F_n-\Phi\Vert_1 \le\frac{E|X_1|^3}{\sigma^3 \sqrt{n}}\qquad
\mbox{for all $n \in\mathbb{N}$,}
\]
corresponding to an $L^1$-Berry--Esseen constant of 1.
\end{abstract}

%

\end{frontmatter}

\section{Introduction}

The classical central limit theorem allows the approximation of the
distribution of sums of ``comparable'' independent real-valued
random variables by
the normal. As this theorem is an asymptotic, it provides no information
as to whether the resulting approximation is useful. For that purpose,
one may turn to the
Berry--Esseen theorem, the most classical version giving supremum norm
bounds between the distribution function
of the normalized sum and that of the standard normal. Various authors
have also considered Berry--Esseen-type bounds
using other metrics and, in particular, bounds in $L^p$. The case
$p=1$, where the value
\[
\Vert F-G\Vert_1=\int_{-\infty}^\infty|F(x)-G(x)|\,dx
\]
is used to measure the distance between distribution functions $F$ and $G$,
is of some particular interest and results using
this metric are known as mean central limit theorems (see, e.g.,
\cite{dedecker,Erickson,HoChen} and \cite{Ibragimov}; the latter
three of these works consider nonindependent summand variables).
One motivation for studying $L^1$-bounds is that, combined with one of
type $L^\infty$, bounds on $L^p$-distance for all $p \in(1,\infty)$
may be
obtained by the inequality
\[
\Vert F-G\Vert_p^p \le\Vert F-G\Vert_\infty^{p-1} \Vert F-G\Vert_1.
\]

For $\sigma\in(0,\infty)$, let $\mathcal{F}_\sigma$ be the collection
of distributions with mean zero, variance $\sigma^2$ and finite absolute
third moment.
We prove the following Berry--Esseen-type result for the mean central
limit theorem.
\begin{theorem}
\label{c1}
For $n \in\mathbb{N}$, let $X_1,\ldots,X_n$ be independent mean zero
random variables with distributions $G_1 \in\mathcal{F}_{\sigma_1},
\ldots, G_n \in\mathcal{F}_{\sigma_n}$ and let
$F_n$ be the distribution of
\[
W =\frac{1}{\sigma}\sum_{i=1}^n X_i\qquad  \mbox{where }
\sigma^2 = \sum_{i=1}^n \sigma_i^2.
\]
Then
\[
\Vert F_n-\Phi\Vert_1 \le\frac{1}{\sigma^3}\sum_{i=1}^n E|X_i|^3.
\]
In particular, when $X_1,\ldots,X_n$ are
identically distributed with distribution \mbox{$G \in\mathcal{F}_\sigma$},
\[
\Vert F_n-\Phi\Vert_1 \le\frac{E|X_1|^3}{\sigma^3 \sqrt{n}}\qquad
\mbox{for all $n \in\mathbb{N}$.}
\]
\end{theorem}

For the case where all variables are identically distributed as $X$
having distribution $G$, letting
%
%
\begin{equation}\label{def-cm}
c_m=\inf\biggl\{C\dvtx \frac{{\sqrt n} \sigma^3 \Vert F_n-\Phi\Vert _1}{E|X|^3}\le C
 \mbox{ for all $G \in\mathcal{F}_1$ and $n \ge m$}\biggr\},
\end{equation}
the second part of Theorem \ref{c1} yields the upper bound $c_1 \le
1$. Regarding lower bounds, we also prove
%
%
\begin{equation}\label{constant}
c_1 \ge\frac{
2 \sqrt{\pi} (2\Phi(1)-1) -
(\sqrt{\pi}+\sqrt{2})+ 2 e^{-1/2}\sqrt{2}
}{\sqrt{\pi}} =0.535377\ldots.
\end{equation}

Clearly, the elements of the sequence $\{c_m\}_{m \ge1}$ are
nonnegative and decreasing in $m$, and so have a limit, say $c_\infty
$. Regarding limiting behavior, Esseen \cite{esseen} showed that
\[
\lim_{n \rightarrow\infty} n^{1/2} \Vert F_n-\Phi\Vert_1 =A(G)
\]
for an explicit constant $A(G)$ depending only on $G$. Zolotarev \cite
{Zolotarev}
provides the representation
%
%
\begin{equation}\label{def-AG}
A(G)=\frac{1}{\sigma\sqrt{2 \pi}}\int_{-1/2}^{1/2} \int_{-\infty
}^\infty\biggl\vert\frac{\omega}{2}(1-x^2)+hu \biggr\vert e^{-x^2/2}\,dx \,du,
\end{equation}
where $\omega=|EX^3|/(3\sigma^2)$ and $h$ is the span of the
distribution $G$ in the case where $G$ is lattice, and is zero otherwise.
Zolotarev obtains
\[
\sup_{G \in\mathcal{F}_\sigma} \frac{\sigma^3 A(G)}{E|X|^3} =
\frac{1}{2},
\]
showing that $c_\infty= 1/2$, hence giving the asymptotic
$L^1$-Berry--Esseen constant value.

Here, the focus is on nonasymptotic constants and, in particular, on
the constant $c_1$ which gives a bound for all $n \in\mathbb{N}$.
Theorem \ref{c1} is shown using Stein's method (see
\cite{Stein72,Stein86}), which uses the characterizing equation (\ref
{Stein-ch}) for the normal, and an associated differential equation to
obtain bounds on the normal approximation. More particularly, we employ
the zero bias transformation, introduced in~\cite{Goldstein-Reinert},
and the evaluation of a Stein functional, as in \cite{Utev}; see, in
particular, Proposition~4.1 there. In \cite{Goldstein-Reinert}, it was
shown that for all $X$ with mean zero and finite nonzero
variance~$\sigma^2$, there exists a unique distribution for a random variable
$X^*$ such that
%
%
\begin{equation}\label{zb-def}
\sigma^2 Ef'(X^*)=E[Xf(X)]
\end{equation}
for all
absolutely continuous functions $f$ for which these expectations
exist. The zero bias transformation, mapping the distribution of
$X$ to that of $X^*$, was motivated by the Stein characterization
of the normal distribution \cite{Stein81}, which states that $Z$ is
normal with
mean zero and variance $\sigma^2$ if and only if
%
%
\begin{equation}
\label{Stein-ch}
\sigma^2 Ef'(Z)=E[Zf(Z)]
\end{equation}
for all absolutely
continuous functions $f$ for which these expectations exist.
Hence, the mean zero normal with variance $\sigma^2$ is the unique
fixed point of the zero bias transformation. How closeness to
normality may be measured by the closeness of a distribution to
its zero bias transform, and related applications, are the topics of
\cite{res,zsm} and \cite{L1bounds}.

As shown in \cite{L1bounds} and \cite{Goldstein-Reinert}, for a
random variable $X$ with $EX=0$ and
$\operatorname{Var}(X)=\sigma^2$,
the distribution of $X^*$ is absolutely continuous
with density and distribution functions given, respectively, by
%
%
\begin{equation}
\label{density-star}
g^*(x)=\sigma^{-2} E[X\mathbf{1}(X>x)]  \quad\mbox{and}\quad
G^*(x)=\sigma^{-2}E[X(X-x)\mathbf{1}(X \le x)].\hspace*{-18pt}
\end{equation}
Theorem \ref{c1} results by showing that the functional
%
%
\begin{equation}\label{def-BG}
B(G)=\frac{2\sigma^2 \Vert G^*-G\Vert_1 }{E|X|^3}
\end{equation}
is bounded by 1 for all $X$ with distribution $G \in\mathcal
{F}_\sigma$.
As in (\ref{def-AG}), one may write out a more ``explicit'' form for
$B(G)$ using (\ref{density-star}) and expressions for
the moments on which $B(G)$ depends, but such expressions appear to be
of little value for the purposes of
proving Theorem \ref{c1}. In turn, the proof here employs convexity
properties of $B(G)$ which
depend on the behavior of the zero bias transformation on mixtures. We
also note that the
functional $B(G)$ has a different character
than $A(G)$; for instance, $A(G)$ is zero for all
nonlattice distributions with vanishing third moment, whereas $B(G)$ is
zero only for mean zero normal
distributions.

Let $\mathcal{L}(X)$ denote the distribution of a random variable $X$.
Since the $L^1$-distance scales, that is, since,
for all $a \in\mathbb{R}$,
%
%
\begin{equation}\label{L1scales}
\Vert \mathcal{L}(aX)-\mathcal{L}(aY)\Vert_1=|a|\Vert\mathcal{L}(X)-\mathcal
{L}(Y)\Vert_1,
\end{equation}
by replacing $\sigma_i^2$ by $\sigma_i^2/\sigma^2$ and
$\Vert G_i^*-G_i\Vert_1$ by $\Vert G_i^*-G_i\Vert_1/\sigma$ in equation (16) of
Theorem 2.1 of \cite{L1bounds}, we obtain the following.
\begin{proposition}
\label{already-shown}
Under the hypotheses of Theorem \ref{c1},
we have
\[
\Vert F_n-\Phi\Vert_1 \le\frac{1}{\sigma^3}\sum_{i=1}^n B(G_i) E|X_i|^3.
\]
\end{proposition}

For $\mathcal{F}$ a collection of nontrivial mean zero distributions
with finite absolute third moments, we let
\[
B(\mathcal{F})=\sup_{G \in\mathcal{F}} B(G).
\]
Clearly, Theorem \ref{c1} follows immediately from Proposition \ref
{already-shown} and the following result.
\begin{lemma} For all $\sigma\in(0,\infty)$,
\label{Gsone}
\[
B(\mathcal{F}_\sigma) =1.
\]
\end{lemma}

The equality to 1 in Lemma \ref{Gsone} improves the upper bound of 3
shown in \cite{L1bounds}. Although our
interest here is in best universal constants, we note that
Proposition~\ref{already-shown} shows that $B(G)$ is
a distribution-specific $L^1$-Berry--Esseen constant, in that
\[
\Vert F_n-\Phi\Vert_1 \le\frac{B(G)E|X_1|^3}{\sigma^3 \sqrt{n}}\qquad
\mbox{for all $n \in\mathbb{N}$},
\]
when $X_1,\ldots,X_n$ are identically distributed according to $G \in
\mathcal{F}_\sigma$.
For instance, $B(G)=1/3$ when $G$ is a mean zero uniform distribution
and $B(G)=1$
when $G$ is a mean zero two-point distribution; see Corollary 2.1 of
\cite{L1bounds}, and Lemmas \ref{lemma-list}
and~\ref{scaling-invariant} below.

We close this section with two preliminaries. The first collects some facts
shown in \cite{L1bounds}
and the second demonstrates that to
prove Lemma \ref{Gsone}, it suffices to consider the class of random
variables $\mathcal{F}_1$.
Then, following Hoeffding \cite{Hoeffding} (see also \cite{Utev}), in
Section \ref{reduction}, we use
a continuity property of $B(G)$ to show that its supremum over
$\mathcal{F}_1$ is attained on
finitely supported distributions. Exploiting a convexity-type property
of the zero bias transformation on mixtures
over distributions having equal variances, we reduce the calculation
further to the calculation of the supremum over~$D_3$,
the collection of all mean zero distributions with variance 1
and supported on
at most three points. As three-point distributions are, in general,
a~mixture of two two-point
distributions with unequal variances, an additional argument is given
in Section \ref{D3}, where a coupling of an $X$ with
distribution $G \in D_3$ to a variable $X^*$ having the $X$ zero bias
distribution is constructed, using the optimal $L^1$-couplings on the
component two-point distributions of which $G$ is the mixture, in order
to obtain $B(G) \le1$ for all $G \in D_3$. The lower bound (\ref
{constant}) on $c_1$
is calculated in Section \ref{lower-bound}.

The following simple formula
will be of some use. For $a \ge0, b>0$ and $l>0$,
we have
%
%
\begin{equation}\label{two-triangles}
\int_0^l  \biggl|(a+b)\frac{u}{l}-a \biggr|\,du = \frac{l}{2}\frac
{a^2+b^2}{a+b}.
\end{equation}
\begin{lemma}
\label{lemma-list}
Let $G$ be the distribution of a nontrivial mean zero random variable
$X$ supported on the two points $x<y$. Then $X^*$ is
uniformly distributed on $[x,y]$,
\[
EX^2=-xy, \qquad  E|X^3|=\frac{-xy(y^2+x^2)}{y-x}
\]
and
\[
\Vert\mathcal{L}(X^*)-\mathcal{L}(X)\Vert_1=\frac{1}{2}\frac{y^2+x^2}{y-x}.
\]
In particular, $B(G)=1$ and
\[
B(\mathcal{F}_1) \ge1.
\]
\end{lemma}
\begin{pf}
Being nontrivial, $G$ has positive variance and,
from (\ref{density-star}), we see that the density $g^*$ of $G^*$ at
$u$, which is
proportional to $E[X\mathbf{1}(X>u)]$, is zero outside $[x,y]$ and
constant within it, so
$G^*(w)=(w-x)/(y-x)$ for $w \in[x,y]$.
That $G$ has mean zero implies that
the support points $x$ and $y$ satisfy $x<0<y$ and that $G$ gives
positive probabilities $y/(y-x)$ and $-x/(y-x)$ to $x$ and $y$,
respectively. The moment
identities are immediate.

Making the change of variable $u=w-x$ and applying (\ref
{two-triangles}) with $a=y/(y-x), b=-x/(y-x)$ and $l=y-x$ yields
\[
\Vert\mathcal{L}(X^*)-\mathcal{L}(X)\Vert_1 = \int_{x}^{y} \biggl|\frac
{w-x}{y-x}-\frac{y}{y-x} \biggr|\,dw
= \frac{1}{2}  \biggl( \frac{y^2+x^2}{y-x}  \biggr),
\]
and (\ref{def-BG}) now gives $B(G)=1$.
\end{pf}
\begin{lemma} \label{scaling-invariant}
Let $G \in\mathcal{F}_\sigma$ for some $\sigma\in(0,\infty)$, let
$X$ have distribution $G$
and, for $a \not= 0$, let $G_a$ denote the distribution of $aX$. Then
$B(G_a)=B(G)$ and, in particular,
\[
B(\mathcal{F}_\sigma)=B(\mathcal{F}_1)  \qquad\mbox{for all $\sigma
\in(0,\infty)$.}
\]
\end{lemma}
\begin{pf}
That $aX^*$ has the same distribution as $(aX)^*$ follows from (\ref
{zb-def}). The identities
$\sigma_{aX}^2=a^2\sigma_X^2$, $E|aX|^3=|a|^3E|X^3|$ and (\ref
{L1scales}) now imply the first
claim. The second claim now follows from
\[
\{B(G)\dvtx G \in\mathcal{F}_\sigma\} = \{B(G)\dvtx G \in\mathcal{F}_1\}.
\]
\upqed\end{pf}

\section{Reduction to three-point distributions}
\label{reduction}

Let $(S,\Sigma)$ be a measurable space
and let $\{m_s\}_{s \in S}$ be a collection of
probability measures on $\mathbb{R}$ such that for each Borel
subset $A \subset\mathbb{R}$, the function from $S$ to
$[0,1]$ given by
\[
s \rightarrow m_s(A)
\]
is measurable. When $\mu$ is a probability measure on
$(S,\Sigma)$, the set function given by
\[
m_\mu(A)=\int_S m_s(A)\mu(d s)
\]
is a probability measure, called the $\mu$ mixture of $\{m_s\}_{s \in
S}$. With
some slight abuse of notation, we let $E_\mu$ and $E_s$
denote expectations with respect to $m_\mu$ and $m_s$, and let
$X_\mu$ and $X_s$ be random variables with distributions
$m_\mu$ and $m_s$, respectively. For instance, for all
functions $f$
which are integrable with respect to $\mu$, we have
\[
E_\mu f(X) =
\int E_s f(X) \mu(d s),
\]
which we also write as
\[
E f(X_\mu) = \int E f(X_s) \mu(d s).
\]
In particular, if $\{m_s\}_{s \in S}$ is a
collection of mean zero
distributions with variances $\sigma_s^2=EX_s^2$ and absolute third
moments
$\gamma_s=E |X_s^3|$, the mixture
distribution $m_\mu$ has variance $\sigma_\mu^2$ and third absolute moment
$\gamma_\mu$ given by
\[
\sigma_\mu^2 = \int_S \sigma_s^2 \,d\mu \quad\mbox{and}\quad
\gamma_\mu= \int_S \gamma_s \,d\mu,
\]
where both may be
infinite. Note that $\sigma_\mu^2 < \infty$ implies that
$\sigma_s^2<\infty$ $\mu$-almost surely and, therefore, that
$m_s^*$, the $m_s$ zero bias distribution, exists
$\mu$-almost surely.

Theorem \ref{mix-zero} shows that the zero bias distribution of a
mixture is a mixture of zero bias distributions, the
mixing measure of which is the original measure weighted by the
variance and rescaled. Define
(arbitrarily, see Remark \ref{arbitrarily}) the zero bias distribution of
$\delta_0$, a point mass at zero, to be $\delta_0$. Write $X =_d Y$ when
$X$ and $Y$ have the same distribution.
\begin{theorem}
\label{mix-zero}
Let $\{m_s, s \in S\}$ be a collection of
mean zero distributions on $\mathbb{R}$ and $\mu$ a probability
measure on $S$ such that the variance $\sigma_\mu^2$ of the
mixture distribution is positive and finite. Then $m_\mu^*$, the
$m_\mu$ zero
bias distribution, exists and is given by the mixture
\[
m_\mu^*=\int m_s^* \,d\nu\qquad \mbox{where }
\frac{d\nu}{d\mu} = \frac{\sigma_s^2}{\sigma_\mu^2}.
\]
In particular, $\nu=\mu$
if and only if $\sigma_s^2$ is a constant $\mu$ a.s.
\end{theorem}
\begin{pf}
The distribution $m_\mu^*$ exists as
$m_\mu$ has mean zero and finite nonzero variance.
Let $X_\mu^*$ have the $m_\mu$ zero bias distribution
and let $Y$ have distribution~$m_\mu^*$. For any absolutely continuous
function $f$ for which the
expectations below exist, we have
\begin{eqnarray*}
\sigma_\mu^2 Ef'(X_\mu^*)&=& EX_\mu f(X_\mu)\\
&=& \int EX_s f(X_s)\,d\mu\\
&=& \int\sigma_s^2 Ef'(X_s^*)\,d\mu\\
&=& \sigma_\mu^2 \int Ef'(X_s^*)\,d\nu\\
&=&\sigma_\mu^2 Ef'(Y).
\end{eqnarray*}
Since $Ef'(X_\mu^*)=Ef'(Y)$ for all such $f$, we conclude that $X_\mu
^*=_d Y$.
\end{pf}
\begin{remark} \label{arbitrarily}
If $m_s = \delta_0$ for any $s \in S$,
then $\sigma_s^2=0$ and, therefore,
\[
\nu\{s \in S\dvtx m_s = \delta_0\}=0.
\]
Hence, the mixture $X_\mu^*$ gives zero weight to all corresponding
$m_s^*$, showing that $(\delta_0)^*$ may be defined arbitrarily.
\end{remark}

We now recall an equivalent form of the $L^1$-distance involving
expectations of Lipschitz functions $L$ on $\mathbb{R}$,
%
%
\begin{eqnarray}\label{lip}
\Vert F-G\Vert_1 = {\sup_{f \in L}}|Ef(X)-Ef(Y)|\nonumber\\[-8pt]\\[-8pt]
\eqntext{\mbox{where }   L=\{
f\dvtx |f(x)-f(y)| \le|x-y|\},}
\end{eqnarray}
$X$ and $Y$ having distributions $F$ and $G$,
respectively. With a slight abuse of notation, we may write $B(X)$
in place of $B(G)$ when $X$ has distribution $G$.
\begin{theorem}
\label{X-is-mixture} If $X_\mu$ is the $\mu$ mixture of a
collection $\{X_s, s \in S\}$ of mean zero, variance 1
random variables satisfying $E|X_\mu^3| <
\infty$, then
%
%
\begin{equation}\label{BxmulsBXalpha}
B(X_\mu) \le\sup_{s \in S} B(X_s).
\end{equation}

If $\mathcal{C}$ is a collection of mean zero, variance 1 random variables
with finite absolute third moments and $\mathcal{D} \subset\mathcal{C}$
such that every distribution in $\mathcal{C}$ can be represented as a
mixture of distributions
in $\mathcal{D}$, then
%
%
\begin{equation}\label{BDldBD}
B(\mathcal{C})=B(\mathcal{D}).
\end{equation}
\end{theorem}
\begin{pf}
Since the variances $\sigma_s^2$ of $X_s$ are
constant, the distribution $X_\mu^*$ is the $\mu$ mixture
of $\{X_s^*, s \in S\}$, by Theorem \ref{mix-zero}. Hence, applying
(\ref{lip}), we have
%
%
\begin{eqnarray}\label{int-alpha}
\Vert\mathcal{L}(X_\mu^*)-\mathcal{L}(X_\mu)\Vert_1 &=& {\sup_{f \in
L}}|Ef(X_\mu^*)-Ef(X_\mu)|\nonumber\\
&=& \sup_{f \in L} \biggl|\int_S Ef(X_s^*)\,d\mu-\int_S
Ef(X_s)\,d\mu \biggr|\nonumber\\[-8pt]\\[-8pt]
&\le& \sup_{f \in L}\int_S  | Ef(X_s^*)- Ef(X_s)
 | \,d\mu\nonumber\\
&\le& \int_S \Vert\mathcal{L}(X_s^*)- \mathcal
{L}(X_s)\Vert_1 \,d\mu.\nonumber
\end{eqnarray}
Noting that $\operatorname{Var}(X_\mu)=\int_S EX_s^2 \,d\mu=1$ and applying
(\ref{int-alpha}),
we find
that
\begin{eqnarray*}
B(X_\mu)&=&\frac{2\Vert\mathcal{L}(X_\mu^*)-\mathcal{L}(X_\mu
)\Vert_1}{E|X_\mu^3|}\\
&\le& \frac{\int_S 2\Vert\mathcal{L}(X_s^*)-\mathcal{L}(X_s)\Vert_1\,d\mu
}{E|X_\mu^3|}\\
&=& \frac{\int_S B(X_s) E|X_s^3| \,d\mu}{E|X_\mu^3|}\\
&\le& \sup_{s \in S}B(X_s) \frac{\int_S E|X_s^3|\,d\mu}{E(X_\mu
^3)}\\
&=& \sup_{s \in S} B(X_s).
\end{eqnarray*}

Regarding (\ref{BDldBD}), clearly, $B(\mathcal{D}) \le B(\mathcal
{C})$ and the reverse
inequality follows from~(\ref{BxmulsBXalpha}).
\end{pf}
%
%
\begin{remark}
Note that no bound of the type provided by Theorem \ref{X-is-mixture}
holds, in general,
when taking mixtures of variables that have unequal
variances. In particular, if $X_s \sim\mathcal{N}(0,\sigma_s^2)$ and
$\sigma_s^2$
is not constant in $s$, then $X_\mu$ is a mixture of normals with
unequal variances, which is not normal.
Hence, in this case, $B(X_\mu) >0$, whereas $B(X_s)=0$ for all $s$.
\end{remark}

To apply Theorem \ref{X-is-mixture} to reduce the computation of
$B(\mathcal{F}_1)$
to finitely supported distributions, we apply the following continuity
property of
the zero bias transformation; see Lemma 5.2 in \cite{lin}. We write
$X_n \Rightarrow X$ for
the convergence of $X_n$ to $X$ in distribution.
\begin{lemma}
\label{zb-cont} Let $X$ and $X_n, n=1,2,\ldots,$ be mean zero
random variables with finite, nonzero variances. If
\[
X_n \dd X  \quad\mbox{and}\quad  \lim_{n \rightarrow
\infty} EX_n^2 = EX^2,
\]
then
\[
X_n^* \dd X^*.
\]
\end{lemma}

For a distribution function $F$, let
%
%
\begin{equation}\label{defFinv}
F^{-1}(w)=\sup\{a\dvtx F(a)<w\}  \qquad\mbox{for all $w \in(0,1)$}.
\end{equation}
If $U$ is uniform on $[0,1]$, then $F^{-1}(U)$ has distribution
function $F$. If $X_n$ and $X$ have
distribution functions $F_n$ and $F$, respectively, and $X_n
\Rightarrow X$, then $F_n^{-1}(U) \rightarrow F^{-1}(U)$ a.s.
(see, e.g., Theorem 2.1 of \cite{durrett}, Chapter 2). For
distribution functions $F$ and $G$,
we have
%
%
\begin{equation}\label{minL1coup}
\Vert F-G\Vert_1 = \inf E|X-Y|,
\end{equation}
where the infimum is over all joint distributions on $X,Y$ which have
marginals $F$ and $G$, respectively,
and the variables $F^{-1}(U)$ and $G^{-1}(U)$ achieve the
minimal $L^1$-coupling, that is,
%
%
\begin{equation}\label{minL1coupU}
\Vert F-G\Vert_1 = E|F^{-1}(U)-G^{-1}(U)|;
\end{equation}
see \cite{Rachev} for details.

With the use of Lemma \ref{zb-cont}, we are able to prove the following
continuity property of the functional $B(X)$.
\begin{lemma}
\label{continuity}
Let $X$ and $X_n$, $n \in\mathbb{N}$, be mean zero
random variables with finite, nonzero absolute
third moments. If
%
%
\begin{equation}\label{12mom-conv}
X_n \dd X,\qquad  \lim_{n \rightarrow
\infty} EX_n^2 = EX^2  \quad\mbox{and}\quad   E|X_n^3| \rightarrow E|X^3|,
\end{equation}
then
\[
B(X_n) \rightarrow B(X) \qquad \mbox{as $n \rightarrow\infty$}.
\]
\end{lemma}
\begin{pf}
By Lemma \ref{zb-cont}, we have $X_n^* \Rightarrow X^*$. Let $U$ be a uniformly
distributed variable and set
\[
(Y,Y_n,Y^*,Y_n^*)=(F_X^{-1}(U),F_{X_n}^{-1}(U),F_{X^*}^{-1}(U),F_{X_n^*}^{-1}(U)),
\]
where $F_W$ denotes the distribution function of $W$.
Then $Y=_d X$, $Y_n=_d X_n$, $Y^*=_d X^*$ and $Y_n^*=_d X_n^*$.
Furthermore, $Y_n \as Y$, $Y_n^* \as
Y^*$ and, by (\ref{minL1coupU}),
\[
\Vert\mathcal{L}(X_n^*)-\mathcal{L}(X_n)\Vert_1=E|Y_n^*-Y_n|  \quad\mbox
{and}\quad  \Vert\mathcal{L}(X^*)-\mathcal{L}(X)\Vert =E|Y^*-Y|.
\]
By (\ref{zb-def}) with $f(x)=x^2\mbox{sgn}(x)$, we find, for $Y$, for
example, that
\[
E|Y^3|=2\operatorname{Var}(Y)E|Y^*|.
\]
Hence, as $n \rightarrow\infty$, we have $EY_n^2 = EX_n^2 \rightarrow
EX^2= EY^2$ and
\[
E|Y_n^*| =\frac{E|Y_n^3|}{2EY_n^2} = \frac
{E|X_n^3|}{2EX_n^2}\quad\rightarrow\quad\frac{E|X^3|}{2EX^2}= \frac
{E|Y^3|}{2EY^2}=E|Y^*|  \qquad\mbox{as $n \rightarrow\infty$}.
\]
Hence, $\{Y_n\}_{n \in\mathbb{N}}$ and $\{Y_n^*\}_{n \in\mathbb
{N}}$ are uniformly integrable, so
$\{Y_n^*-Y_n\}_{n \in\mathbb{N}}$ is uniformly integrable. As
$Y_n^*-Y_n \as Y^*-Y$ as $n\rightarrow\infty$,
we have
%
%
\begin{eqnarray}\label{cont-of-L1}
{\lim_{n \rightarrow\infty}} \Vert\mathcal{L}(X_n^*)-\mathcal
{L}(X_n)\Vert_1 &=& \lim_{n \rightarrow\infty} E|Y_n^*-Y_n| = E|Y^*-Y| \nonumber\\[-8pt]\\[-8pt]
&=&
\Vert\mathcal{L}(X^*)-\mathcal{L}(X)\Vert.\nonumber
\end{eqnarray}
Combining (\ref{cont-of-L1}) with the convergence of the variances and
the absolute third moments, as provided by (\ref{12mom-conv}),
the proof is complete.
\end{pf}

Lemmas \ref{cX-xY} and \ref{m->3} borrow much from Theorem 2.1 of
\cite{Hoeffding}, the latter
lemma indeed being implicit. However, the results of \cite{Hoeffding}
cannot be applied directly as $B(G)$ is not expressed as the
expectation of $K(X)$ for some
$K$ when $\mathcal{L}(X)=G$. For $m \ge2$, let $D_m$ denote the
collection of all mean
zero, variance 1 distributions which are supported on at most $m$ points.
\begin{lemma}
\label{cX-xY}
\[
B(\mathcal{F}_1)= B\biggl(\bigcup_{m \ge3}D_m\biggr).
\]
\end{lemma}
\begin{pf}
Letting $\mathcal{M}$ be the collection of distributions in $\mathcal
{F}_1$ which have compact support, we first show that
%
%
\begin{equation}\label{bf-le-bm}
B(\mathcal{F}_1) \le B(\mathcal{M}).
\end{equation}
Let $\mathcal{L}(X) \in\mathcal{F}_1$ be given and, for $n \in
\mathbb{N}$, set
$Y_n=X\mathbf{1}_{|X| \le n}$. Clearly, $Y_n \dd X$. As
$E|X^3|<\infty$ and $|Y_n^p| \le|X^p|$ for all $p \ge0$, by the
dominated convergence theorem, we have
%
%
\begin{eqnarray}\label{YntoX}
EY_n &\rightarrow& EX=0,\nonumber\\[-8pt]\\[-8pt]
EY_n^2 &\rightarrow& EX^2=1  \quad\mbox
{and}\quad   E|Y_n^3| \rightarrow E|X^3|  \qquad\mbox{as $n \rightarrow
\infty$}.\nonumber
\end{eqnarray}
Letting
%
%
\begin{equation}\label{def-Xn}
X_n=Y_n-EY_n,
\end{equation}
we have $X_n \Rightarrow X$, by Slutsky's theorem, so, in view of (\ref
{YntoX}), the hypotheses of
Lemma \ref{continuity} are satisfied, yielding
\[
B(X_n) \rightarrow B(X)  \qquad\mbox{as $n \rightarrow\infty$, with
$\{X_n\}_{n \in\mathbb{N}} \subset\mathcal{M}$},
\]
showing (\ref{bf-le-bm}).

Now, consider $\mathcal{L}(X) \in\mathcal{M}$ so that $|X|\le M$
a.s. for some $M>0$.
For each $n \in\mathbb{N}$, let
\[
Y_n = \sum_{k \in\mathbb{Z}} \frac{k}{2^n}\mathbf{1}\biggl(\frac{k-1}{2^n}
< X \le\frac{k}{2^n}\biggr).
\]
Since $|X| \le M$ a.s., each $Y_n$ is supported on finitely many points
and uniformly bounded.
Clearly, $Y_n \rightarrow X$ a.s. and (\ref{YntoX}) holds by the
bounded convergence theorem. Now, defining $X_n$ by (\ref{def-Xn}),
the hypotheses of
Lemma \ref{continuity} are satisfied, yielding
\[
B(X_n) \rightarrow B(X)  \qquad\mbox{as $n \rightarrow\infty$, with
}\{X_n\}_{n \in\mathbb{N}} \subset\bigcup_{m \ge3}D_m,
\]
showing $B(\mathcal{M}) \le B(\bigcup_{m \ge3}D_m)$.
Combining this inequality with (\ref{bf-le-bm}) yields $B(\mathcal
{F}_1) \le B(\bigcup_{m \ge3}D_m)$ and therefore the lemma, the
reverse inequality being obvious.
\end{pf}
\begin{lemma}
\label{m->3} Every distribution in $\bigcup_{m \ge3} D_m$
can be expressed as a finite mixture of $D_3$ distributions.
\end{lemma}
\begin{pf}
The lemma is trivially true for $m=3$, so consider $m>3$ and
assume that the lemma holds for all integers from $3$ to
$m-1$.

The distribution of any $X \in D_m$ is determined by the
supporting values $a_1 < \cdots<a_m$ and a vector of probabilities
$\mathbf{p}=(p_1,\ldots,p_m)'$.
If any of the components of $\mathbf{p}$ are zero, then $X \in D_k$ for
$k<m$ and the induction would
be finished, so we assume that all components of $\mathbf{p}$ are strictly
positive. As $X \in D_m$, the vector $\mathbf{p}$ must
satisfy
\[
A\mathbf{p}=\mathbf{c}\qquad  \mbox{where }   A=
\left[
\matrix{
1 & 1 & \cdots& 1\cr
a_1 & a_2 & \cdots& a_m\cr
a_1^2 & a_2^2 & \cdots& a_m^2}
\right]
\mbox{ and }  \mathbf{c}=\left[
\matrix{
1\cr
0\cr
1}
\right].
\]
Since $A \in\mathbb{R}^{3 \times m}$ with $m
> 3$, $\mathcal{N}(A) \not= \{0\}$, that is, there
exists $\mathbf{v} \not= 0$ with
%
%
\begin{equation}\label{Av=0} A\mathbf{v}=0.
\end{equation}
Since
$\mathbf{v} \not=0$ and the equation specified by the first row of
$A$ is $\sum_i v_i=0$, the vector $\mathbf{v}$ contains
both positive and negative numbers. Since the vector $\mathbf{p}$ has
strictly positive components, the numbers $t_1$ and $t_2$ given by
\[
t_1=\inf\Bigl\{t>0\dvtx \min_i(p_i+tv_i) \ge0 \Bigr\}
\quad\mbox{and}\quad
t_2=\inf\Bigl\{t>0\dvtx \min_i(p_i-tv_i) \ge0 \Bigr\}
\]
are both strictly positive. Note that
\[
\mathbf{p}_1=\mathbf{p}+t_1\mathbf{v}  \quad\mbox{and}\quad  \mathbf{p}_2=\mathbf{p}-t_2\mathbf{v}
\]
satisfy
\[
A\mathbf{p}_1=A(\mathbf{p}+t_1\mathbf{v})=A\mathbf{p}=\mathbf{c} = A\mathbf{p}=A(\mathbf{p}-t_2\mathbf{v})=A\mathbf{p}_2,
\]
by (\ref{Av=0}), so that $\mathbf{p}_1$ and $\mathbf{p}_2$ are probability vectors
since their components are nonnegative and sum to one. Additionally, the
corresponding distributions have mean zero and
variance 1, and in each of these two vectors, at
least one component has been set to zero. Hence, we may express the
$m$-point probability vector $\mathbf{p}$ as the mixture
\[
\mathbf{p}=\frac{t_2}{t_1+t_2}\mathbf{p}_1 + \frac{t_1}{t_1+t_2}\mathbf{p}_2
\]
of probability vectors on at most $m-1$ support points, thus showing
$X$ to be the mixture of two distributions in $D_{m-1}$, completing
the induction.
\end{pf}

The following theorem is an immediate consequence of Theorem \ref
{X-is-mixture} and Lemmas \ref{cX-xY}
and \ref{m->3}.
\begin{theorem} \label{F1=D3}
\[
B(\mathcal{F}_1) = B(D_3).
\]
\end{theorem}

Hence, we now restrict our attention to $D_3$.

\section{Bound for $D_3$ distributions}
\label{D3}

Lemma \ref{scaling-invariant} and Theorem \ref{F1=D3} imply that
$B(\mathcal{F}_\sigma)=B(\mathcal{F}_1)=B(D_3)$.
Hence, Lemma \ref{Gsone} follows from Theorem \ref{D3=1} below, which
shows that $B(D_3)=1$.
We prove Theorem \ref{D3=1} with the help of the following result.
\begin{lemma}
\label{QPslePPs}
Let $x<y<0<z$, and let $m_1$ and $m_0$ be the unique mean zero
distributions with support $\{x,z\}$ and $\{y,z\}$, respectively, that is,
\[
m_1(\{w\}) = \cases{
\dfrac{z}{z-x}, &\quad if $w=x$,\vspace*{2pt}\cr
\dfrac{-x}{z-x}, &\quad if $w=z$,\vspace*{2pt}\cr
0, &\quad otherwise,}
\quad\mbox{and}\quad
m_0(\{w\}) = \cases{
\dfrac{z}{z-y}, &\quad if $w=y$,\vspace*{2pt}\cr
\dfrac{-y}{z-y}, &\quad if $w=z$,\vspace*{2pt}\cr
0, &\quad otherwise.}
\]
Then
%
%
\begin{equation}\label{desired}
\Vert m_1^*-m_0\Vert_1 \le\Vert m_1^*-m_1\Vert_1.
\end{equation}
\end{lemma}
\begin{pf}
Let $F_1,F_0$ and $F_1^*$ denote the distribution functions of
$m_1,m_0$ and~$m_1^*$, respectively.
By Lemma \ref{lemma-list}, $m_1^*$ is uniform over $[x,z]$. There are
two cases, depending on the relative magnitudes of
$F_1^*(y)=(y-x)/(z-x)$ and $F_0(y)=z/(z-y)$.

We first consider the case
%
%
\begin{equation}\label{cases}
F_1^*(y) \le F_0(y)  \quad\mbox{or, equivalently,}\quad   y(x+z) \le y^2+z^2.
\end{equation}
By Lemma \ref{lemma-list},
%
%
\begin{eqnarray}\label{m1-m1s}
\Vert m_1^*-m_1\Vert_1&=&\frac{z^2+x^2}{2(z-x)}\nonumber\\
&=&\frac
{(z^2+x^2)(z-y)^2}{2(z-x)(z-y)^2}
\nonumber\\[-8pt]\\[-8pt]
&=& \frac{z^4-2yz^3+y^2z^2+x^2z^2-2x^2yz+x^2y^2}{2(z-x)(z-y)^2}\nonumber\\
&=& \frac{(z^4-2yz^3+x^2z^2-2x^2yz)+y^2z^2+x^2y^2}{2(z-x)(z-y)^2}.\nonumber
\end{eqnarray}
Letting $J_1=[x,y)$ and $J_2=[y,z]$, we have
\[
\Vert m_1^*-m_0\Vert_1 = I_1 + I_2\qquad  \mbox{where }   I_i=\int
_{J_i}|F_1^*(w)-F_0(w)| \,dw  \qquad\mbox{for $i\in\{1,2\}$}.
\]
Since $F_1^*(w) \ge0= F_0(w)$ for all $w \in J_1$,
%
%
\begin{equation}\label{Isub1}
I_1 = \int_x^y  \biggl( \frac{w-x}{z-x}  \biggr) \,dw = \frac
{1}{2}\frac{(y-x)^2}{z-x}= \frac{(y-x)^2(z-y)^2}{2(z-x)(z-y)^2}.
\end{equation}
Recalling that $F_1^*(y) \le F_0(y)$, applying (\ref{two-triangles}) with
$a=\frac{z}{z-y}-\frac{y-x}{z-x}$, $b=-\frac{y}{z-y}$ and $l=z-y$, after
the change of variable $u=w-y$, yields
%
%
\begin{eqnarray}\label{Isub2}\qquad
I_2 &=&\int_y^z  \biggl|\frac{w-x}{z-x}- \frac{z}{z-y} \biggr|\,dw \nonumber\\
&=& \biggl( \frac{z-y}{2}  \biggr) \frac{({z}/({z-y})-
({y-x})/({z-x}))^2+({y}/({z-y}))^2}{1-({y-x})/({z-x})}\nonumber\\
&=& \frac{1}{2}(z-x) \biggl( \biggl(\frac{z}{z-y}-\frac{y-x}{z-x}\biggr)^2+\biggl(\frac
{y}{z-y}\biggr)^2  \biggr)\\
&=& \frac{(z(z-x)-(y-x)(z-y))^2+(y(z-x))^2}{2(z-x)(z-y)^2}\nonumber\\
&=& \frac{(y^2+z^2)(z-x)^2 - 2z(z-x)(y-x)(z-y)+ (y-x)^2(z-y)^2}{2(z-x)(z-y)^2}.\nonumber
\end{eqnarray}
Adding (\ref{Isub1}) to (\ref{Isub2}) yields
\begin{eqnarray*}
&&
\Vert m_1^*-m_0\Vert_1\\
&&\qquad=\frac{(y^2+z^2)(z-x)^2 - 2z(z-x)(y-x)(z-y)+
2(y-x)^2(z-y)^2}{2(z-x)(z-y)^2}\\
%
%
&&\qquad=\frac
{(z^4-2yz^3+x^2z^2-2x^2yz)+5y^2z^2+3x^2y^2-4xy^3}{2(z-x)(z-y)^2}\\
&&\qquad\quad{} +
\frac{4xy^2z-4xyz^2+2y^4-4y^3z}{2(z-x)(z-y)^2}.
\end{eqnarray*}
Now, subtracting from (\ref{m1-m1s}) and simplifying by noting that
the terms
inside the parentheses in the numerators of these two expressions are
equal, we find
that
%
%
\begin{eqnarray}\label{difference}
&&
\Vert m_1^*-m_1\Vert_1 - \Vert m_1^*-m_0\Vert_1 \nonumber\\
&&\qquad= \frac
{-4y^2z^2-2x^2y^2+4xy^3-4xy^2z+4xyz^2-2y^4+4y^3z}{2(z-x)(z-y)^2}\\
&&\qquad=\frac{-y(y-x)(y^2+2z^2-y(x+2z))}{(z-x)(z-y)^2}.\nonumber
\end{eqnarray}
The denominator in (\ref{difference}) is positive, as is $-y$ and $y-x$.
For the remaining term, (\ref{cases}) yields
\[
y^2+2z^2-y(x+2z) = y^2+2z^2-yz-y(x+z)\ge z(z-y) >0.
\]
Hence, (\ref{difference}) is positive, thus proving (\ref{desired})
when $F_1^*(y) \le F_0(y)$.

When $F_1^*(y) > F_0(y)$, we have $F_1^*(w) \ge F_0(w)$ for all $w \in[x,z)$
as $F_0(w)$ is zero in $[x,y)$ and equals $F_0(y)$ in $[y,z)$, and
$F_1^*(w)$ is increasing over $[y,z)$.
Hence,
\begin{eqnarray*}
\Vert m_1^*-m_0\Vert_1&=& \int_x^z |F_1^*(w)-F_0(w)| \,dw
= \int_x^z \bigl(F_1^*(w) - F_0(w)\bigr) \,dw \\
&=& \int_x^z \frac{w-x}{z-x} \,dw - \int_y^z \frac{z}{z-y}\,dw
= \frac{1}{2}\frac{(z-x)^2}{z-x} - z
= \frac{1}{2}(z-x)-z\\
&=&- \frac{x+z}{2}.
\end{eqnarray*}
Now, since $(x+z)(x-z) = x^2-z^2 \le z^2+x^2$ and $z-x>0$, using Lemma
\ref{lemma-list}, we obtain
\[
\Vert m_1^*-m_0\Vert_1=- \frac{x+z}{2} \le\frac{z^2+x^2}{2(z-x)}= \Vert m_0^*-m_0\Vert_1,
\]
thus proving inequality (\ref{desired}) when $F_1^*(y) > F_0(y)$ and,
therefore, proving the lemma.
\end{pf}
\begin{theorem} \label{D3=1}
\[
B(D_3)=1.
\]
\end{theorem}
\begin{pf}
Lemma \ref{lemma-list} shows that $B(X) =1$ if $X$ is supported on two
points, so
$B(D_3) \ge1$ and it only remains to consider $X$ positively supported
on three points.
We first prove that
%
%
\begin{eqnarray}\label{nonzero-support}
B(X) \le1  \hspace*{200pt}\nonumber\\[-8pt]\\[-8pt]
\eqntext{\mbox{when $X \in D_3$ is positively supported on the
nonzero points $x,y,z$.}}
\end{eqnarray}
$EX=0$ implies that $x<0<z$. After proving (\ref{nonzero-support}), we
treat the remaining
case, where $y=0$, by a continuity argument.

Let $X$ be supported on $x<y<z$ with $y \not=0$. Lemma \ref
{scaling-invariant} with $a=-1$ implies that $B(-X)=B(X)$, so we may
assume, without loss of generality, that $x < y < 0 <z$. Let $m_1$ and
$m_0$ be the unique mean zero distributions supported on $\{x,z\}$ and
$\{y,z\}$, respectively, and let $\mathcal{L}(X_1)=m_1$ and $\mathcal
{L}(X_0)=m_0$. As, in general, every mean zero distribution having no
atom at zero can be represented as a mixture of mean zero
two-point distributions
(as in the Skorokhod representation, see \cite{durrett}), letting
%
%
\begin{equation}\label{Xalpha}
\mathcal{L}(X_\alpha)=\alpha m_1 + (1-\alpha)m_0,
\end{equation}
we have $\mathcal{L}(X)=\mathcal{L}(X_\alpha)$ for some $\alpha\in
[0,1]$; in fact, for the given $X$, one may verify that
$P(X=x)/P(X_1=x) \in(0,1)$
and that (\ref{Xalpha}) holds when $\alpha$ assumes this value. Therefore,
to prove (\ref{nonzero-support}), it suffices to show that
%
%
\begin{equation}
\label{bound-for-all-alpha}
B(X_\alpha) \le1  \qquad\mbox{for all $\alpha\in[0,1]$.}
\end{equation}

By Lemma \ref{lemma-list},
%
%
\begin{equation}\label{sigmas-10}
EX_1^2=-zx  \quad\mbox{and}\quad   EX_0^2=-zy,
\end{equation}
and, by (\ref{Xalpha}), the variance of $X_\alpha$ is given by
%
%
\begin{eqnarray}\label{sec-mom-alpha}
EX_\alpha^2 &=& \alpha EX_1^2 + (1-\alpha) EX_0^2 = - \bigl( \alpha zx
+ (1-\alpha) zy  \bigr)\nonumber\\[-8pt]\\[-8pt]
&=&-z\bigl(\alpha x + (1-\alpha) y\bigr).\nonumber
\end{eqnarray}
Applying Theorem \ref{mix-zero} with $S=\{0,1\}$ and $\mu$ being the
probability measure
putting mass $\alpha$ and $1-\alpha$ on the points $1$ and $0$,
respectively, in view
of (\ref{sigmas-10}) and (\ref{sec-mom-alpha}), $m_\alpha^*$, the
$X_\alpha$ zero
bias distribution, is given by the mixture
%
%
\begin{equation}\label{def-beta}
m_\alpha^*=\beta m_1^* + (1-\beta)m_0^*
\qquad\mbox{where }
\beta= \frac{\alpha x}{\alpha x + (1-\alpha)y}.
\end{equation}
Since $x<y<0$, we have
\[
\frac{\beta}{1-\beta} = \frac{\alpha}{1-\alpha}\frac{x}{y}>
\frac{\alpha}{1-\alpha}  \quad\mbox{and, therefore,}\quad  \beta
>\alpha.
\]
Let $F_1,F_0,F_1^*$ and $F_0^*$ denote the distribution functions of
$m_1,m_0,m_1^*$ and $m_0^*$, respectively.
Let $U$ be a standard uniform
variable and, with the inverse functions below given by (\ref
{defFinv}), set
\[
(Y_1,Y_0,Y_1^*,Y_0^*)=(F_1^{-1}(U), F_0^{-1}(U),
(F_1^*)^{-1}(U),(F_0^*)^{-1}(U)).
\]
Then $Y_i=_d X_i$, $Y_i^*=_d X_i^*$ for $i \in\{1,2\}$ and, by (\ref
{minL1coupU}), all pairs of the variables $Y_1,Y_0,Y_1^*,Y_0^*$
achieve the $L^1$-distance between their respective distributions.
Now, let $(Y_\alpha,Y_\alpha^*)$ be defined on the same space with
joint distribution given by the mixture
\[
\mathcal{L}(Y_\alpha,Y_\alpha^*)= \alpha\mathcal{L}(Y_1,Y_1^*) +
(1-\beta)\mathcal{L}(Y_0,Y_0^*)+(\beta-\alpha)\mathcal{L}(Y_0,Y_1^*).
\]
Then $(Y_\alpha,Y_\alpha^*)$ has marginals $Y_\alpha=_d X_\alpha$
and $Y_\alpha^*=_d Y_\alpha^*$, hence, by (\ref{minL1coup}),
%
%
\begin{eqnarray}
\label{coupling-bound}
\Vert m_\alpha^*-m_\alpha\Vert_1 &\le& \alpha\Vert m_1^*-m_1\Vert_1 + (1-\beta
)\Vert m_0^*-m_0\Vert_1\nonumber\\[-8pt]\\[-8pt]
&&{} + (\beta-\alpha)\Vert m_1^*-m_0\Vert_1.\nonumber
\end{eqnarray}

Lemma \ref{lemma-list} shows that $G(X_i)=1$, that is,
$E|X_i^3|=2EX_i^2\Vert m_i^*-m_i\Vert_1$ for $i=1,2$, so (\ref{Xalpha}) yields
\[
E|X_\alpha^3|=2\bigl( \alpha EX_1^2\Vert m_1^*-m_1\Vert_1 + (1-\alpha)
EX_0^2\Vert m_0^*-m_0\Vert_1\bigr)
\]
and, by (\ref{sigmas-10}), (\ref{sec-mom-alpha}) and (\ref
{def-beta}), we now find that
%
%
\begin{eqnarray}
\label{fit}
\frac{E|X_\alpha^3|}{2EX_\alpha^2}&=&\frac{\alpha x\Vert m_1^*-m_1\Vert_1
+(1-\alpha) y \Vert m_0^*-m_0\Vert_1}{\alpha x +
(1-\alpha)y}\nonumber\\[-8pt]\\[-8pt]
&=&\beta\Vert m_1^*-m_1\Vert_1 +(1-\beta)\Vert
m_0^*-m_0\Vert_1.\nonumber
\end{eqnarray}
Lemma \ref{QPslePPs} shows that the right-hand side and, therefore,
the left-hand side of (\ref{coupling-bound}) are bounded by (\ref
{fit}), that is,
that $B(X_\alpha) = 2EX_\alpha^2 \Vert m_\alpha^*-m_\alpha
\Vert_1/E|X_\alpha^3|\le1$, completing
the proof of (\ref{bound-for-all-alpha}) and hence of (\ref{nonzero-support}).

Finally, we consider the case where the mean
zero random variable $X$ is positively supported on $\{x,0,z\}$ with
$x<0<z$ and
$P(X=0)=q \in(0,1)$.
For $n \in\mathbb{N}$, let
\[
Y_n=X+n^{-1}\mathbf{1}(X=0)  \quad\mbox{and}\quad   X_n=Y_n-EY_n.
\]
As $n \rightarrow\infty$, we see that $Y_n \as X$ and
$EY_n =q/n \rightarrow0$ so that
$X_n \as X$, and the bounded convergence theorem shows
that $\{X_n\}_{n \in\mathbb{N}}$ satisfies the hypothesis of Lemma
\ref{continuity}. Hence,
$B(X_n) \rightarrow B(X)$ as $n \rightarrow\infty$. For all $n \in
\mathbb{N}$ such that $1/n<z$,
the distribution of $X_n$ is positively supported on the three
distinct, nonzero
points $x-q/n<(1-q)/n<z-q/n$, so, by (\ref{nonzero-support}),
$B(X_n) \le1$ for all such $n$. Therefore,
the limit $B(X)$ is also bounded by 1.
\end{pf}

\section{Lower bound}
\label{lower-bound}

By (\ref{def-cm}), with $m=1$ and $\mathcal{L}(X)=G \in\mathcal{F}_1$,
\[
\Vert F_n-\Phi\Vert_1 \le\frac{c_1E|X^3|}{\sqrt{n}}  \qquad\mbox{for all
$n \in\mathbb{N}$}
\]
and, in particular, for $n=1$,
%
%
\begin{equation}\label{c1lower}
c_1 \ge\frac{\Vert F_1-\Phi\Vert_1}{E|X^3|}=\frac{\Vert G-\Phi\Vert_1}{E|X^3|}.
\end{equation}

Motivated by Theorem \ref{F1=D3}, that two-point distributions achieve
the suprema of $B(G)$, for $p \in(0,1)$, let
\[
X=\frac{\xi-p}{\sqrt{pq}},
\]
where $\xi$ is a Bernoulli variable with $P(\xi=1)=p=1-P(\xi=0)$.
The distribution function $G_p$ of $X$ is given by
\[
G_p(x)= \cases{
0, &\quad for $\displaystyle x \le-\sqrt{\frac{p}{q}}$,\vspace*{1pt}\cr
q, &\quad for $\displaystyle -\sqrt{\frac{p}{q}} < x \le\sqrt{\frac{q}{p}}$,\vspace*{1pt}\cr
1, &\quad for $\displaystyle \sqrt{\frac{q}{p}}<x$,}
\]
and, therefore, the $L^1$-distance between $G_p$ and the standard
normal is given by
\[
\Vert G_p-\Phi\Vert_1=\int_{-\infty}^{-\sqrt{{p/q}}} \Phi(x)\,dx +
\int_{-\sqrt{{p/q}}}^{\sqrt{{q/p}}}
|\Phi(x)-q|\,dx+\int_{\sqrt{{q/p}}}^\infty|\Phi(x)-1|\,dx.
\]
As $G_p \in\mathcal{F}_1$ for all $p \in(0,1)$ and $E|X^3|=
(p^2+q^2)/\sqrt{pq}$,
letting
\[
\psi(p)=\frac{\sqrt{pq}}{p^2+q^2}\Vert G_p-\Phi\Vert_1  \qquad\mbox{for $p
\in(0,1)$},
\]
inequality (\ref{c1lower}) gives $c_1 \ge\psi(p)$ for all $p \in
(0,1)$ and
$\psi(1/2)$ yields (\ref{constant}).

\section{Remarks}
This article was submitted on November 18th, 2008. In the article \cite
{Tyu}, submitted on June 8th, 2009,
Ilya Tyurin independently proved Theorem~\ref{c1}, also by applying
the zero bias method. The current article was posted on arXiv on June
28, 2009; article \cite{Tyu} was posted on December 3rd, 2009.
In \cite{Tyu}, Theorem \ref{c1} is used to prove the upper bound
$0.4785$ on the $L^\infty$-Berry--Esseen constant.

\section*{Acknowledgment}
The author would like to sincerely thank Sergey Utev for helpful
suggestions.


%
\printaddresses

\end{document}